\documentclass[12pt]{article}

\begin{document}
\title{On the critical dissipative quasi-geostrophic equation}

\author{Peter Constantin \thanks{Partially supported by NSF DMS-9802611}\\Department of Mathematics\\The University of Chicago\\\\ Diego Cordoba\\Department of Mathematics\\ The University of Chicago\\\\ 
 Jiahong Wu\thanks{Partially supported by NSF-DMS 9971926.} 
\\Department of Mathematics\\Oklahoma State University\\
 }

\date{}
\maketitle

\newtheorem{thm}{Theorem}[section]
\newtheorem{cor}[thm]{Corollary}
\newtheorem{prop}[thm]{Proposition}
\newtheorem{define}[thm]{Definition}
\newtheorem{rem}[thm]{Remark}
\newtheorem{example}[thm]{Example}
\newtheorem{lemma}[thm]{Lemma}
\def\theequation{\thesection.\arabic{equation}}

\noindent{\bf Abstract}: 
The 2D quasi-geostrophic (QG) equation 
is a two dimensional model of the 3D incompressible Euler equations.
When dissipation is included in the model then
solutions always exist if the dissipation's wave number dependence is 
super-linear. Below this critical power the
dissipation appears to be insufficient. For instance, it is not known if 
the critical dissipative QG equation has  
global smooth solutions for arbitrary large initial data.
In this paper we prove existence and uniqueness of global classical 
solutions of the critical dissipative QG equation for initial data 
that have small $L^\infty$ norm. The importance of an $L^{\infty}$ smallness
condition is due to the fact that $L^{\infty}$ is a conserved norm for the
non-dissipative QG equation and is non-increasing on all solutions of
the dissipative QG., irrespective of size.

\vspace{.3in}
\noindent AMS classification scheme numbers: 86A10, 35Q35, 76U05

\section{Introduction}
\setcounter{equation}{0}
\label{sec:1}

Do singularities develop in finite time in smooth 
solutions of unforced, incompressible 3D fluid  equations? This challenging question remains yet unanswered.  Lower dimensional model equations have been 
proposed and studied (\cite{CLM},\cite{Ma},\cite{CMT},\cite{Co}) in an 
attempt to develop mathematical insight in  this problem. The 2D 
quasi-geostrophic (QG) equation is one of these models.
The dissipative QG equation is
\begin{equation}\label{qgse}
\frac{\partial \theta}{\partial t}+u\cdot\nabla\theta +
\kappa (-\Delta)^{\alpha}\theta=0,
\end{equation}
where $\alpha\in[0,1]$, $\kappa>0$ is 
the dissipative coefficient, and the 2D velocity field
$u=(u_1, u_2)$ is determined from $\theta$ by a stream function 
 $\psi$ via the auxiliary relations
\begin{equation}\label{utheta}
(u_1,u_2)=\left (-\frac{\partial \psi}{\partial x_2}, \frac{
\partial \psi}{\partial x_1}\right), \qquad 
{(-\Delta)^{\frac{1}{2}}}\psi = \theta.
\end{equation}
Equation (\ref{qgse}) with $\alpha=1/2$ is the critical dissipative QG. 
Criticality means that the dissipation balances nonlinearity when one takes
into account the conservation laws.

In addition to its intrinsic mathematical interest, the equation (\ref{qgse}) 
is relevant in the context of general quasi-geostrophic 
models of atmospheric and ocean fluid flow \cite{Pe}.

We are concerned here with  global existence results for solutions of 
the initial-value
problem (IVP) for equation (\ref{qgse}) wherein 
\begin{equation}\label{init}
\theta(x,0) =\theta_0(x) 
\end{equation}
is specified. We consider periodic boundary conditions with period box 
$\Omega=[0,2\pi]^2$. Without loss of generality 
we may restrict the discussion to $\theta$ that obey
for all time  $(2\pi)^{-2}\int_{\Omega} \theta dx = 0$.

The issue of global existence for equation (\ref{qgse}) is non-trivial
(\cite{CW1},\cite{Wu}). If no smallness condition is imposed on the 
initial data, then the issue of global existence for arbitrary data is open.  
In this paper we show that if the $L^\infty$-norm of the initial data
is small then (\ref{qgse}) possesses a 
global solution in the critical case $\alpha=1/2$. The QG equations
(dissipative or not) have global weak solutions for arbitrary $L^2$
initial data (\cite{R}, see Appendix B). The $L^{\infty}$ norm condition is 
significant in view of the fact that the QG equations have a maximum principle
(\cite{R}, see Appendix A) that ensures that for all time and initial data
the $L^{\infty}$ norm is non-increasing in time. 
If the initial data is smooth enough then the solution of the 
critical dissipative QG is unique, smooth and decays in time. These 
results and their proofs are presented in section \ref{sec:2}.  
When $\alpha>1/2$, the smallness assumption on the data is not needed
for global existence of smooth solutions for equation 
(\ref{qgse}). This can be proved using the same ideas as 
for the critical case. The  theory of global existence and regularity 
in the sub-critical $(\alpha>1/2)$ case 
is thus in a satisfactory state. More details for the sub-critical case 
can be found in 
\cite{Wu}.

We establish now some of the notation. 
The Fourier transform of $f$ is $\widehat{f}$ 
$$
\widehat{f}(k) =\frac1{(2\pi)^2 }\int_{\Omega} f(x) e^{-ik\cdot x} dx.
$$
$\Lambda$ is used to denote the operator 
$(-\Delta)^\frac{1}{2}$, defined at the Fourier level by  
$$
\widehat{\Lambda f}(k) = |k| \widehat{f}(k).
$$
The relation in (\ref{utheta}) can be identified as 
\begin{equation}\label{utheta1}
u=\left(-\partial_{x_2} \Lambda^{-1} \theta,\, \partial_{x_1} 
\Lambda^{-1} \theta\right)  = (-{R}_2 \theta,\, {R}_1\theta)
\quad \mbox{or}\quad \widehat{u}(j) = i\left (\widehat{j}\right )^{\perp}
\widehat{\theta}(j), 
\end{equation}
where $i = \sqrt{-1}$, $\widehat{j} = \frac{j}{|j|}$ for $j\in{\mathbf Z}^2\setminus\{0\}$, $(j_1,j_2)^\perp = (-j_2,j_1)$,  
and ${R}_1$ and ${R}_2$ are Riesz transforms (\cite{St}).
The spaces $H^s$ are the familiar Sobolev spaces of functions having $s$ derivatives in $L^2$.

\vspace{.2in}
\section{Global existence}
\setcounter{equation}{0}
\label{sec:2}

The initial value problem for the critical dissipative QG equation is 
\begin{equation}\label{dqg}
\left\{
\begin{array}{ll}
\theta_t  + u\cdot \nabla \theta 
+ \kappa \Lambda \theta=0,\quad 
& (x,t)\in \Omega\times [0,\infty),\\\cr
u=(u_1,u_2) = (-{R}_2 \theta,\, {R}_1\theta),\quad
& (x,t)\in \Omega\times [0,\infty),\\\cr
\theta(x,0) =\theta_0(x),\quad
& x\in \Omega,\cr
\end{array}
\right.
\end{equation}
where $\kappa>0$ is a constant. 

\vspace{.1in}
In this section we 
assume that the $L^\infty$ norm of the initial data 
$\theta_0$ is small. We establish that
the IVP (\ref{dqg}) has  a global bounded solution
in $H^1$. If the initial data is smoother
($H^2$) then the solution's norm in $H^2$ is non-increasing in time. The
solution becomes real analytic at positive time and and decays exponentially.  
For the sub-critical case, no smallness assumption is necessary and these results hold for arbitrary data.

We start with an apriori estimate.
\begin{thm}\label{apr}
There exists a constant $c_{\infty}$ such that
for any $\theta_0\in H^2\cap{\mathcal C}^3$ with 
\begin{equation}\label{small}
\|\theta_0\|_{L^\infty} \le c_{\infty} \kappa,
\end{equation}
the classical solution $\theta$ of the IVP (\ref{dqg}) 
satisfies 
\begin{equation}\label{bdd}
\|\theta(\cdot,t)\|_{H^2} \le \|\theta_0\|_{H^2}
\end{equation}
for all $t\ge 0$. 
\end{thm}
{\it Proof.}\quad 
Multiplying the first equation in (\ref{dqg}) by $\Delta^2\theta$ and 
integrating by parts, we obtain
$$
\frac12\frac{d}{dt} \int |\Delta\theta|^2\, dx  + \kappa \int |(-\Delta)^\frac54
\theta|^2 dx  = -\int \Delta^2 \theta\, (u\cdot\nabla \theta) dx. 
$$
Further integration by parts gives 
$$
\int \Delta^2 \theta\, (u\cdot\nabla \theta) dx
= 2 \int \nabla u\cdot (\nabla (\nabla\theta)) \Delta \theta + 
\int (\Delta u \cdot \nabla \theta) \Delta\theta
$$ 
By H\"{o}lder's inequality, 
$$
\left|\int \Delta^2 \theta\, (u\cdot\nabla \theta) dx \right| 
\le C \left[\|\nabla u\|_{L^3} \|\Delta \theta\|_{L^3}^2  + \|\Delta u\|_{L^3}
\|\nabla \theta\|_{L^3} \|\Delta \theta\|_{L^3}\right]
$$
The Riesz transforms are bounded in $L^p$ spaces, so
$$
\|\Delta u\|_{L^3} \le C \|\Delta \theta\|_{L^3},
\qquad \|\nabla u \|_{L^3} \le C \|\nabla \theta\|_{L^3}.
$$ 
The Gagliardo-Nirenberg inequalities
$$
\|\nabla \theta\|_{L^3} \le C \|\theta\|_{L^\infty}^\frac79 
\|(-\Delta)^\frac54\theta\|^\frac29_{L^2}, \qquad
\|\Delta \theta\|_{L^3} \le C \|\theta\|_{L^\infty}^\frac19 
\|(-\Delta)^\frac54\theta\|^\frac89_{L^2},
$$
follow from classical ones (\cite{gn}) by complex interpolation (see also{\cite{kp}}). Using them we obtain 
$$
\left|\int (-\Delta)^2 \theta (u\cdot\nabla \theta) dx \right|
\le C \|\theta\|_{L^\infty} \|(-\Delta)^\frac54\theta\|^2_{L^2}.
$$
Collecting the above estimates, we have 
$$
\frac12 \frac{d}{dt}\int |\Delta\theta|^2\, dx  + \kappa \int |(-\Delta)^\frac54
\theta|^2 dx \le C_{\infty} \|\theta\|_{L^\infty} \|(-\Delta)^\frac54\theta\|^2_{L^2}.
$$
It was proved in (\cite{R}) (see the Appendix A) that $\theta$ 
satisfies the maximum principle 
$$
\|\theta(\cdot,t)\|_{L^\infty} \leq \|\theta_0\|_{L^\infty} 
\quad\mbox{for all $t\ge 0$}. 
$$
Taking $c_{\infty} = (C_{\infty})^{-1}$ the bound (\ref{bdd}) then follows 
from the smallness condition (\ref{small}).
This completes the proof of the theorem.

\vspace{.15in}
\begin{thm}
There exists a constant $c_{\infty}$ 
(the same as in Theorem \ref{apr} above) so that for any
$\theta_0\in H^2$ with $\|\theta_0\|_{L^\infty}\le c_{\infty}\kappa$ 
the IVP (\ref{dqg}) has 
a unique global solution $\theta$ satisfying 
$$
\|\theta(\cdot,t)\|_{H^2} \le \|\theta_0\|_{H^2}
$$
for any $t\ge 0$.
\end{thm}
{\it Proof}.\quad Let $\theta$ be the unique local solution on $[0,T_0]$ with $T_0$
depending on $\|\theta_0\|_{H^2}$ only (standard techniques can be applied to
show that a unique local solution 
exists and depends continuously on initial data in $H^2$). 
By Theorem \ref{apr}, $\theta$ satisfies 
$$
\|\theta(\cdot,t)\|_{H^2} \le \|\theta_0\|_{H^2}
$$
for any $t\in [0,T_0]$. Therefore the local solution can be extended uniquely to 
$[0,2T_0]$ and the global solution is obtained by repeating
this procedure.

\begin{thm} Assume that the initial data $\theta_0\in H^2$ satsifies the
bound $\|\theta_0 \|_{L^\infty} <c_{\infty}\kappa$. Then the solution of the IVP 
(\ref{dqg}) decays exponentially 
$$
\|\theta\|^2_{H^2} \le \exp {(-ct)}\|\theta_0\|^2_{H^2}
$$
for all $t\ge 0$. Here $c = 2(\kappa - c_{\infty}^{-1}\|\theta_0\|_{L^{\infty}})$,
and $c_{\infty}$ is the same as in the previous theorems.
\end{thm}
The proof is a trivial consequence of the obvious Poincare inequality 
$$
\int_{\Omega} |(-\Delta)^\frac 54
\theta|^2 dx \ge \int_{\Omega} |(-\Delta)
\theta|^2 dx.
$$

\begin{thm} Assume that the initial data $\theta_0\in H^2$ satsifies the
bound $\|\theta_0 \|_{L^\infty} <c_{\infty}\kappa $. Then there exists a time 
$t_0 >0 $ such that the solution of the IVP (\ref{dqg}) is real analytic 
for $t\ge t_0$.
More precisely, there exists an extension of $\theta$,
$\theta (z, t)$ that is an analytic function in the 
time-expanding strip $\Sigma_t\subset {\mathbf C}^2$,
$$
\Sigma_t = \{z = x+iy;\,\, x\in \Omega,\,\, |y| <\frac{1}{2}\kappa (t-t_0)\}
$$
and obeys the inequality
$$
|\theta (z,t) | \le \frac{\kappa}{2}
$$
uniformly for $z\in \Sigma_t$.
\end{thm}

{\it Proof.} Let us consider $t_0$ to be the first time when
$$
Y(t) = \sum_{j\in{\mathbf Z}^2\setminus \{0\}}|\widehat{\theta}(j,t)|
$$ 
becomes smaller than $\kappa/4$:
$$
Y(t_0)\le \kappa/4.
$$
In view of the preceding theorem and the elementary inequality
$$
Y(t) \le C\|\theta (t)\|_{H^2}
$$
the existence of $t_0$ is guaranteed. Consider the function
$$
y(t) = \sum_{j\in{\mathbf Z}^2\setminus \{0\}}|\widehat{\theta}(j,t)|
\exp{\left \{\frac{(t-t_0)\kappa |j|}{2}\right \}}
$$
and the function
$$
z(t) = \sum_{j\in{\mathbf Z}^2 \setminus \{0\}}|j||\widehat{\theta}(j,t)|\exp{\left \{\frac{(t-t_0)\kappa |j|}{2}\right \}}
$$
Formally
$$
\frac{dy}{dt} + z(t)(\frac{\kappa}{2} - y(t)) \le 0
$$
holds. This implies that the function $y(t)$ is non-increasing for $t\ge t_0$.
The more rigorous proof requires one to take only a finite sum, and introduce an
artificial power to avoid differentiating the modulus at zero:
$$
y_{n,\epsilon}(t) = \sum_{j\in{\mathbf Z}^2\setminus \{0\}, \, |j|\le n}
|\widehat{\theta}(j,t)|^{1+\epsilon}\exp{\left \{\frac{(t-t_0)\kappa |j|}{2}\right \}}
$$
This is now a differentiable function in time. One differentiates, and obtains
$$
\frac{d}{dt}y_{n,\epsilon}(t) +\frac{\kappa}{2}z_{n,\epsilon}(t) \le \gamma_{\epsilon} 
y(t)z(t)
$$
where 
$$
z_{n,\epsilon}(t) =  \sum_{j\in{\mathbf Z}^2\setminus \{0\}, \, |j|\le n}
|j| |\widehat{\theta}(j,t)|^{1+\epsilon}\exp{\left \{\frac{(t-t_0)\kappa |j|}{2}\right \}}, 
$$
and $\gamma_{\epsilon} = (1+\epsilon)\|\theta_0\|_{L^2}^{\epsilon}$, and
thus  $\lim_{\epsilon\to 0}\gamma_{\epsilon} = 1$. One integrates from
$t=t_0$ to $t$ and passes to the limit $n\to\infty$, $\epsilon\to 0$.
One obtains then
$$
y(t) + \frac{\kappa}{2}\int_{t_0}^tz(s)ds \le \frac{\kappa}{4} + \int_{t_0}^ty(s)z(s)ds
$$
Then, because $y(t_0) \le \frac{\kappa}{4}$ it follows that the set
$\{t\ge t_0|\,\, y(s)\le \frac{\kappa}{2},\,\forall  s, t_0\le s\le t\}$  equals $[t_0, \infty)$.
The completely rigorous proof requires a regularization of the equation
so that $z(t)$ is guaranteed to be finite. After $t=t_0$ one may use Galerkin approximations for this purpose. 
The conclusion is that, for $t\ge t_0$ one has
$$
\sum_{j\in{\mathbf Z}\setminus \{0\}}|\widehat{\theta}(j,t)| \exp{\left \{\frac{(t-t_0)\kappa |j|}{2}\right \}}\le \frac{\kappa}{2}
$$
The analytic extension is
$$
\theta(z,t) = \sum_{j\in{\mathbf Z}^2\setminus \{0\}}e^{i j\cdot z}\widehat{\theta}(j,t)
$$
and the uniform convergence and bound in $\Sigma_t$ follow. This completes the
proof. The idea of using time dependent exponential weights was introduced
in (\cite{gev}).

\vspace{.15in}
Initial data in  $H^1$ are  sufficient for a  
global existence results to hold.
\begin{thm}\label{h1}
There exists a constant $d_{\infty}$ such that
for any 
$\theta_0\in H^1$ and $\|\theta_0\|_{L^\infty} \le d_{\infty}\kappa$
there exists a weak solution of the QG equation satisfying
$$
\|\theta(\cdot,t)\|_{H^1} \le \|\theta_0\|_{H^1}
$$
for any $t\ge 0$.
\end{thm}
{\it Proof}.\quad 
It was proved in (\cite {R}) that the IVP (\ref{dqg}) with 
$\theta_0\in L^2$ has a global weak 
solution $\theta(\cdot,t)\in L^2$ (see Appendix B). The $L^2$ weak solutions are constructed using a Galerkin approximation. The $H^1$ weak solutions
can be constructed by solving approximate 
equations
\begin{equation}
\partial_t\theta + u_{\delta}\cdot\nabla \theta + \kappa \Lambda \theta = 0, \label{app}
\end{equation}
where
$$
u_{\delta} = k_{\delta}*u= k_{\delta}*\left (R^{\perp}\theta\right )
$$
and $k_{\delta }$ the periodic Poisson kernel in 2D given at the Fourier level by
$$
\widehat{k_{\delta}}(\xi) = e^{-\delta |\xi|}
$$
$\xi\in {\mathbf Z}^2$.
The approximations have gobal smooth solutions for positive time, uniform
bounds in $L^{\infty}(dt; L^2(dx))$ and converge weakly to solutions
of the QG equation. In addition, and in contrast with Galerkin approximations, these 
approximations  have monotonic non-increasing $L^p$ norms of $u$.
Multiplying the first equation in (\ref{app}) by $\Delta \theta$
and integrating by parts, we obtain 
$$
\frac12 \frac{d}{dt}\int |\nabla \theta|^2 dx + \kappa \int |\Lambda^\frac32\theta|^2 dx 
\le \int |(\nabla \theta)\cdot\nabla u_{\delta}\cdot (\nabla \theta)| dx.
$$ 
By a similar argument as in the proof of Theorem \ref{apr}, 
the term on the right hand side can be bounded as follows. 
$$
\int |(\nabla \theta)\cdot\nabla u_{\delta}\cdot (\nabla \theta)| dx 
\le C \|\nabla \theta\|^2_{L^3} \,\|\nabla u_{\delta}\|_{L^3}
\le C \|\nabla \theta\|^3_{L^3} 
\le C \|\theta\|_{L^\infty}\|\Lambda^\frac32\theta\|_{L^2}^2.
$$
In addition to the boundedness of Riesz transforms in $L^p$
spaces and an appropriate Gagliardo-Nirenberg inequality, we use here
the fact that convolution with the Poisson kernel does not 
increase $L^p$ norms.  Therefore,
$$
\frac12 \frac{d}{dt}\int |\nabla \theta|^2 dx + \kappa \int |\Lambda^\frac32\theta|^2 dx
\le D_{\infty}\|\theta\|_{L^\infty}\|\Lambda^\frac32\theta\|_{L^2}^2.
$$
Using $d_{\infty} = (D_{\infty})^{-1}$ and the maximum principle
we deduce from
$\|\theta_0\|_{L^\infty} \le d_{\infty}$ and the 
inequality above that  
\begin{equation}\label{time}
\|\theta(\cdot,t)\|_{H^1} \le \|\theta_0\|_{H^1} \quad\mbox{and}\quad
\int_0^t \int |\Lambda^\frac32\theta|^2(x,\tau) dx\,d\tau < \infty
\end{equation}
for any $t>0$. The solutions obtained thus are relatively strong and the equation holds in time integral form in $L^2$. 
\vspace{.15in}

\vspace{.15in}
For the sub-critical case $\alpha>\frac12$ the maximum principle allows one to get a sublinear bound of the nonlinearity in terms of the dissipation and,
consequently, the global existence result
holds without any smallness assumption on $\|\theta\|_{L^\infty}$.
\begin{thm}
Let $\alpha>\frac12$ and $\theta_0\in H^2$. Then there exists a unique 
global solution $\theta$ solving the IVP (\ref{dqg}). The solution is real analytic for positive time and decays exponentially to zero.
\end{thm}

\vspace{.2in}

\section{Appendix A}

In \cite{R} it was shown that a  solution to (\ref{dqg})
 for $\frac{1}{2}\le \alpha\le 1$ satisfies the following 
maximum principle 
$$
\|\theta(\cdot,t)\|_{L^p} \leq \|\theta_0\|_{L^p} \quad\mbox{for $1<p\leq \infty
$} \quad\mbox{for all $t\ge 0$}.
$$

Below we give a short description of the proof in the case 
$\alpha = \frac{1}{2}$.

Define $\tilde{\theta} = k_s * \theta$ where $\hat{k}_s(\xi) = e^{-s|\xi|}$. 
The Poisson kernel $k_s$ is positive and has 
integral equal to one. $\tilde{\theta}$ satisfies
$$
\frac{d \tilde{\theta}}{ds} + \Lambda\tilde{\theta} = 0
$$
and
$$
\frac{d \|\tilde{\theta}\|_{L^p}^{p}}{ds} + p \int |\tilde{\theta}|^{p-2}\tilde{
\theta} \Lambda\tilde{\theta} dx = 0.
$$
Integrating with respect to s  we obtain
$$
 p \int_{s_1}^{s_2} \int |\tilde{\theta}|^{p-2}\tilde{\theta} \Lambda\tilde
{\theta} dx ds = \|\tilde{\theta}\|_{L^p}^{p}(s_1) - \|\tilde{\theta}\|_{L^p}^{p}(s_2
) 
$$
Using the  properties of the Poisson Kernel $k_s$

$$
\|k_{s_1} * \theta\|_{L^p}^{p} - \|k_{s_2 - s_1} *(k_{s_1} * \theta)\|_{L^p}^{p}
\geq 0.
$$
therefore
$$
\int |\theta |^{p-2}\theta\Lambda\theta dx = lim_{s\rightarrow 0^+} \int |\tilde
{\theta}|^{p-2}\tilde{\theta} \Lambda\tilde{\theta} dx \geq 0.
$$
Hence, returning to the evolution equation for $\theta$
$$
\frac{d \|\theta\|_{L^p}^{p}}{dt} = - p \int |\theta||^{p-2}\theta\Lambda
\theta dx \leq 0.
$$

\section{Appendix B}
It was proved in (\cite{R}) that the QG equations, dissipative or not, have 
global weak solutions in $L^2$. We present here a brief description of the main 
reason for this fact.
We consider QG as an infinite system of ODEs:
\begin{equation}
\frac{d}{dt}{\widehat{\theta}}(l,t)  + 
\kappa |l|^{2\alpha}{\widehat{\theta}}(l,t)  = b_l(\theta, \theta)\label{feq}
\end{equation}
with
\begin{equation}
b_l(\theta, \theta) = \sum_{j+k = l}\frac{1}{|j|}\left (j^{\perp}\cdot k \right ){\widehat{\theta}}(j,t)  {\widehat{\theta}}(k,t) 
\label{bl}
\end{equation}
Note that because of the perpendicularity
$$
\left (j^{\perp}\cdot k \right ) = \left (j^{\perp}\cdot l\right) = - \left (k^{\perp}\cdot l\right )
$$
and because of symmetry considerations in the sum we can write
\begin{equation}
b_l(\theta, \theta) =  \sum_{j+k = l}\gamma^{l}_{j,k}{\widehat{\theta}}(j,t)  {\widehat{\theta}}(k,t) 
\label{blo}
\end{equation}
with 
$$
\gamma^{l}_{j,k} = \frac{1}{2}\left (j^{\perp}\cdot l\right)\frac{|k|-|j|}{|j||k|}.
$$
Note the inequality 
$$
\left | \gamma^{l}_{j,k} \right | \le \frac{|l|^2}{2(\max{|j|; \, |k|})}
$$ 
Consider now the weak norm
$$
\|\theta\|_{w} = \sup_{j\in {\mathbf Z}^2\setminus \{0\}}|\widehat{\theta }(j)|
$$
The nonlinearity can be written as
$$
B(\theta,\theta )(x) = \sum_{l\in{\mathbf Z}^2\setminus \{0\}}b_l(\theta,\theta)
e^{i l\cdot x}
$$
The main observation is that the nonlinearity has  a weak continuity property.
Let $\theta_1$ and $\theta _2$ be in $L^2$. There exists a constant $C$
so that

$$
\|\Lambda ^{-2}(B(\theta _1, \theta_1) - B(\theta_2, \theta_2))\|_w \le
$$
$$
C\|\theta_1 - \theta _2\|_w\left (1 + \log {(1+ \|\theta_1 - \theta _2\|_w^{-1})}\right )
\left (\|\theta_1\|_{L^2} + \|\theta _2\|_{L^2}\right )
$$
This is the key ingredient in the existence of weak solutions. The method of proof
employs a Galerkin approximation. The weak continuity
guarantees the fact that the weak limit of the approximations solves the
weak  formulation of the equation.

\noindent{\bf Acknowledgment}  We thank A. Kiselev for helpful discussions.

\end{document}